\documentclass[11pt]{article}
\usepackage[utf8]{inputenc}
\usepackage{graphicx}
\usepackage{amsfonts}
\usepackage[width=17cm, left=1cm]{geometry}
\usepackage{amssymb}
\usepackage{amsmath}
\usepackage{amsthm}

\newtheorem{definition}{Definition}
\newtheorem{proposition}{Proposition}
\newtheorem{theorem}{Theorem}

\newtheorem{lemma}{Lemma}
\newtheorem{corollary}{Corollary}
\newtheorem{problem}{Problem}
\title{The asymptotics of the partition of the cube into Weyl simplices, and an encoding of a Bernoulli scheme\thanks{Partially supported by the RFBR grant 17-01-00433.}}

\begin{document}

\author{A.~M.~Vershik\thanks{St.~Petersburg Department of Steklov Institute of Mathematics, St.~Petersburg State University, and Institute for Information Transmission Problems.}}
\date{02.02.2019}
\maketitle
\begin{abstract}
We suggest a combinatorial method of encoding continuous symbolic dynamical systems. A~continuous phase space, the infinite-dimensional cube, turns into the path space of a tree, and the shift is mapped to a transformation which was called a ``transfer.'' The central problem is that of distinguishability: does the encoding separate almost all points of the space? The main result says that the partition of the cube into Weyl simplices satisfies this property.
\end{abstract}

\tableofcontents

\section{Introduction}
In the theory of dynamical systems and, more generally, measure theory, the following question often proves fruitful: find a cover of a given measure space by a space endowed with a product measure (i.e., a product of independent variables) for which a given transformation (or a group of transformations) is a homomorphic image of the shift (respectively, are homomorphic images of the shifts) in the covering space. If such a cover (homomorphism) exists, and even is an isomorphism, then we obtain important information on the original object.

On the other hand, the inverse problem is not less important: how one can economically encode a~sequence of  independent identically distributed  continuous random variables (i.e., a Bernoulli scheme) using finite codes? That is,  can one replace a continuous scheme by a locally finite one and how can this be done?

In both cases, it is important to find a method of economical encoding of a continuous scheme and Bernoulli shift, or even a more general system.

We suggest a general method of \emph{combinatorial encoding} of a Bernoulli scheme, and consider a~simplest nontrivial example of such an encoding, using the partition of the cube into Weyl simplices. This example is related to a simplest tree (of permutations) and produces an isomorphism between the classical Bernoulli scheme with a continuous set of states and a new type of a measure-preserving transformation  called \emph{transfer}. This is a nonstationary Markov shift acting in the Cantor-like space of paths in a graded graph (in the simplest case, in the path space of a tree).

In this context, objects of ergodic theory become related to the combinatorial theory of graded graphs and, consequently, to combinatorics and representation theory.

The main problem arising here is the \emph{distinguishability problem}: is the encoding faithful, i.e., does it separate almost all points of the Bernoulli scheme? In other words, is the encoded system isomorphic to the original one or is it only a homomorphic image of this system? In the example with Weyl simplices, this problem can be stated very simply: can one recover a realization of a Bernoulli scheme with state space~$[0,1]$ from all pairwise inequalities between its coordinates? Quite paradoxically, the answer (see Section~3) is positive:  this can be done with probability~$1$. Actually, this question is related to the theory of equidistributed sequences, but the extensive literature on the subject seems to contain no mention of this fact.

A much more complicated example is related to the Young graph and the RSK correspondence. It was initiated by the old paper~\cite{KV} and was resolved, also positively, in the recent papers~\cite{RS,Sn}. In this case, the encoding uses $Q$-tableaux of the RSK correspondence. The proofs in~\cite{RS,Sn}  are based on a thorough analysis of the theorem on the limit shape of Young diagrams and the study of  the Sch\"utzenberger transformation.

We will return to this example in another paper, where we will apply a general method of resolving the distinguishability problem, which consists in proving certain laws of large numbers ``along'' realizations. These laws can be quite complicated. For example, our proof is based on a new theorem on the limit shape for $P$-tableaux; however, the approach itself, which is mentioned in the present paper too, is universal.

It is worth noting that our notion of transfer (see Sections~4 and~5) is a far generalization of the Sch\"utzenberger transformation  (jeu de taquin) and seems to be important for the general theory of graphs and transformations with an invariant measure. Describing all invariant measures for a transfer and the study of its properties is a new interesting area of the theory of dynamical systems and the combinatorics of graphs.

The main technique here is a combinatorial method of studying increasing invariant sequences of finite measurable partitions of Lebesgue spaces, which  are, according to V.~A.~Rokhlin, separable complete measure spaces, or, in other words, spaces isomorphic mod~0 to an interval with the Lebesgue measure (if the original measure is continuous). In the author's opinion, the most important problems in measure theory and its applications to different fields of mathematics are related to the geometry and combinatorics of $\sigma$-algebras (= measurable partitions) and their sequences. Here we are interested in infinite increasing sequences; the combinatorics underlying the theory of such sequences is the study of properties of infinite trees and graded graphs. Another class of sequences of $\sigma$-algebras is that of decreasing sequences, or filtrations (see~\cite{UMN}); the corresponding theories are closely related, but strongly different.

In Section~2, we define a combinatorial encoding and state the main problems. In particular, we define a frame, i.e., a tree endowed with a translation which is a combinatorial invariant of an increasing sequence of partitions. We state and discuss the distinguishability problem, and also  discuss numerical invariants of exhausting sequences, for which the distinguishability problem has a positive answer. The main Section~3  introduces an encoding of a Bernoulli scheme by Weyl simplices (=~intersections of Weyl chambers with the unit cube); it is used to establish an isomorphism between a~continuous Bernoulli shift and the transfer of a triangular compactum (= the path space of the tree of i-permutations). Section~4 contains various problems related to this example and its generalizations. In particular, we state a problem on different compactifications of the infinite symmetric group, one  being the compactum of virtual permutations and another one being the main example of this paper. In Section~5, we give a general definition of transfer for a graded graph and describe problems related to this notion.

\section{Classical and combinatorial encoding of transformations}

\subsection{Classical encoding and generators}

We begin by recalling the method of encoding endomorphisms and automorphisms or, more generally, arbitrary actions of groups and semigroups with an invariant measure used in symbolic dynamics, when one defines an action of the group by shifts in the space of functions on the group with a~shift-invariant measure.

For simplicity, we consider a Bernoulli endomorphism (= one-sided Bernoulli shift)~$S$, which is a~transformation of the space $\prod_n I\equiv I^{\infty}$, where $I=[0,1]$, with an invariant product measure~$m^{\infty}$, where $m$ is the Lebesgue measure on  $[0,1]$. The usual method of encoding endomorphisms (and automorphisms) in ergodic theory and information theory is to choose an $I$-valued measurable function ${f: I^{\infty}\rightarrow I}$  on the space of trajectories $I^{\infty}$ of the process (``symbolic space'') and study the family of all its shifts $\{f(S^{-n}\cdot)\}_{n=1}^{\infty}$. If these shifts (regarded as functions on~$I^{\infty}$) separate the points of $I^{\infty}$, i.e., the product of the partitions of $I^{\infty}$ into the preimages of points corresponding to the shifts of~$f$ is the partition into singletons, then the partition into the level sets of the original function (and the function itself) is called  a {\it generator}. In this case, we obtain a new isomorphic model of the shift in the same space $I^{\infty}$, but, in general, with another invariant measure different from $m^{\infty}$. To check whether a~given function is a generator is a difficult and instructive problem even for a Bernoulli scheme, and more so in the general case. But since the product of partitions is a partition invariant under the shift, we always have a well-defined quotient endomorphism, which may or may not be isomorphic to the original endomorphism. This classical method of encoding endo- and automorphisms is studied in many papers both in ergodic theory and information theory.

\subsection{Combinatorial encoding}

We suggest to encode an endomorphism (in particular, a Bernoulli shift) in another way: instead of classical codes (functions or partitions), we will construct a {\it shift-invariant increasing sequence of finite measurable partitions on a measure space (in particular, on the infinite-dimensional cube $(I^{\infty}, m^{\infty})$  with a product measure)}.

\smallskip

In more detail, let $I=[0,1]$ and consider the   infinite-dimensional cube $I^{\infty}=\prod_n I=\{\{\xi_n\}_{n\ge1}\}$ with the Lebesgue measure $m^{\infty}=\prod_n m$ where $m$ is the Lebesgue measure on the interval $[0,1]$. By $S$ we denote the one-sided shift, i.e., the endomorphism $S(\{\xi_n\}_{n\geq 1})=\{\xi_{n+1}\}_{n\geq 1}$ of the space $I^{\infty}$ with the invariant measure $m^{\infty}$.

\smallskip

We consider arbitrary increasing  sequences of finite cylinder partitions $\eta_n$, $n=1,2, \dots$, of the  infinite-dimensional cube $I^{\infty}$.

A  {\it finite partition} consists of finitely many (measurable) sets of positive measure, called its elements. A sequence of finite partitions $\{\eta_n\}_n$ is
 {\it increasing} (notation: $\eta_n\prec \eta_{n+1}$) if for every~$n$ every element~$C\in\eta_{n+1}$ is a subset of some element $D\in \eta_n$, and every element $D\in \eta_n$
is the union of all elements of~$\eta_{n+1}$ belonging to it. To exclude degenerate cases, we will assume that every element of~$\eta_n$ contains at least two elements of~$\eta_{n+1}$.

 The  {\it invariance of a sequence of partitions} $\{\eta_n\}_n$  under the shift~$S$ means that the images~$Sx,Sx'$ of (almost) any two points  $x,x'\in I^{\infty}$ belonging to the same element of~$\eta_n$ belong to the same element of $\eta_{n-1}$, for $n=2,3 \dots$.

\smallskip

Besides, we require $\{\eta_n\}_n$ to be cylinder partitions, with $\eta_n$ being a partition into cylinders whose bases are subsets of the finite-dimensional cube~$I^n$, which is the projection of the infinite-dimensional cube to the first $n$ coordinates. The more general case where $\eta_n$ is a cylinder partition with respect to the cube $I^{k_n}$, where $k_n$ is an increasing sequence with $k_n\geq n$, is no different from this one.

Hence, describing a sequence of partitions of the type under consideration reduces to describing a sequence of coherent partitions of  finite-dimensional cubes, which we will denote by the same symbols~$\eta_n$.

\smallskip

We will assume that $\eta_1$ is the trivial partition (whose unique element is a set of full measure); its base is the trivial partition of the first cube, i.e., the interval $I^1=I$.

Assume that we have already constructed an increasing chain of partitions
 ${\eta_1 \prec \eta_2 \prec \dots \prec \eta_n}$ invariant under finitely many shifts in which the elements of~$\eta_k$  are cylinder sets with bases in  ${I^k (\leftarrow I^{\infty})}$, for $k=1,2,\dots, n$. Consider the cube $I^{n+1}$. For the base of the next partition $\eta_{n+1}$, we can take an {\it arbitrary partition of the cube $I^{n+1}$ that is a refinement of the product  $\eta'_n\vee\eta''_n$} of two partitions defined as follows: take two projections $\pi_1,\pi_2: I^{n+1}\rightarrow  I^n$ of the cube
 $I^{n+1}$ to the cube $I^n$, the first one along the axes
$\{1,2, \dots, n\}$, and the second one along the axes $\{2,3, \dots, n+1\}$; then  $\eta'_n$ and $\eta''_n$ are the preimages of the partition $\eta_n$ of $I^n$ under these projections: $$\eta'_n=\pi_1^{-1}(\eta_n),\qquad\eta''_n=\pi_2^{-1}(\eta_n).$$

Continue this process; the fact that the resulting sequence of partitions is increasing and invariant immediately follows from construction. Thus, we have obtained an infinite increasing shift-invariant sequence of finite cylinder partitions. It is not difficult to see that this procedure allows one to construct an arbitrary monotone invariant sequence of finite partitions. One can say that with this method of encoding, the continuity of the state space ``escapes to infinity.'' This phenomenon is worth a more detailed analysis.

Two sequences $\{\eta_n\}_n$ and $\{\eta'_n\}_n$ of the type under consideration are said to be metrically isomorphic if there exists an invertible measurable transformation $T:I^{\infty}\rightarrow I^{\infty}$ preserving the measure $m^{\infty}$ and sending one sequence to the other one.

Note that the described procedure can be applied to an arbitrary Lebesgue space and a measure-preserving transformation of this space defined in symbolic form, i.e., with a generator fixed beforehand.

\subsection{The frame as a combinatorial invariant of an encoding}

With an increasing invariant sequence of finite cylinder partitions
 $\{\eta_n\}_n$ of the space
  $(I^{\infty},m^{\infty})$ we associate a most important combinatorial object: an infinite tree
 $F=F(\{\eta_n\}_{n=1}^{\infty})$ with an additional structure introduced below. The vertices of the $n$th level of~$F$ correspond to the elements of the partition~$\eta_n$; let $C \in \eta_n$ be one of these elements; then it is joined by an edge with the vertex corresponding to the element $D \in \eta_{n-1}$ of the previous partition that contains $C$; we assume that this edge is directed from $D$ to $C$. Thus, we have defined an infinite, locally finite tree corresponding to a~sequence of partitions. Such a tree can be defined for every increasing sequence of finite partitions of every Lebesgue space.

But we have also the following map defined on the vertices of~$F$: since the sequence of partitions~$\{\eta_n\}_n$  is shift-invariant, to every element $C\in \eta_n$ for $n>1$ there corresponds a unique element  $E \in \eta_{n-1}$ different from~$D$ that is the image of all points of~$C$ under the shift  $S:I^{\infty}\rightarrow I^{\infty}$; take the edge connecting the vertices corresponding to the elements $C$ and $E$ and direct it from $E$ to $C$. Thus, our tree is endowed with a bijection from the set of all its vertices to the set of all vertices except the first one, and this bijection preserves the partial order, i.e., sends a pair of vertices that constitute an edge to a similar pair at the previous level. We will call this map a \textbf{translation} and denote by~$\omega$.  Using $\omega$, we define a map from the set of all infinite paths $\{t_i\}_1^{\infty}$ in~$F$ to itself; namely, given a path~$\{t_i\}_1^{\infty}$, the vertices of its image $\{t'_i\}_1^{\infty}$ are the translations of its vertices: $ t'_i=\omega(t_{i+1})$, $i=1,2, \dots$.
This map will be called a \textbf{transfer} on the path space of the tree, and a tree for which a~transfer is defined will be called a \textbf{tree with a transfer}.

Finally, recall that every element of every partition $\eta_n$ has a measure, a positive number from the interval $(0,1)$, and the sum of these numbers over each level of the tree is equal to~$1$. It is more convenient to fix the conditional measure on the elements of the quotient partition $\eta_n/{\eta_{n+1}}$, i.e., fix a~probability vector for each element $C\in \eta_n$.

  \begin{definition}
The frame of a combinatorial encoding of the space $(I^{\infty}, m^{\infty})$, i.e., of an increasing invariant sequence of finite cylinder partitions $\{\eta_n\}_n$, is the tree with a transfer $F(\{\eta_n\}_n)$ defined above endowed with a coherent system of probability vectors on its levels. In the case most interesting for our purposes, the conditional measures are uniform and determined by numbers
$r_n(C)$.
  \end{definition}

The frame, being a graded tree with a transfer endowed with a system of measures, is a combinatorial (finite) invariant of an increasing sequence of partitions  $\{\eta_n\}$ of the space $(I^{\infty},m^{\infty})$. On the other hand, every graded tree with a transfer and a system of measures can be realized as the frame of an increasing sequence of invariant finite measurable partitions of a Lebesgue space with a~measure-preserving transformation  (which  can be different from $(I^{\infty}, m^{\infty}, S)$).

We will say that two sequences of partitions (or two encodings) are combinatorially isomorphic if their frames are isomorphic as graded trees with a transfer and a system of measures.

Clearly, two metrically isomorphic sequences of partitions are combinatorially isomorphic; however, the converse is not true, since the behavior of the sequences at infinity can be different (see the distinguishability problem below).

For the tree of Weyl simplices of type~$A$, the frame and the corresponding transfer are considered in Section~3.

\subsection{Distinguishability problem}

Consider an infinite increasing shift-invariant sequence of finite cylinder partitions $\eta_n$ of the space $I^{\infty}$. Recall that the limit of an increasing sequence of finite partitions $\eta_n$, $ n=1,2,\dots$, is their product, i.e., the measurable partition $\eta\equiv \bigvee_n \eta_n$ whose elements are all nonempty intersections of sequences of elements of~$\eta_n$.

The fundamental question is whether the limiting partition, i.e., the product of partitions, coincides $\mod 0$ with the partition into singletons, i.e., whether it separates almost all, with respect to the measure
$m^{\infty}$, points of the space $I^{\infty}$. If the answer is positive, this means that our encoding loses no information. In this case, the sequence $\{\eta_n\}$ will be called {\it  exhausting}. It is more conventional to say that such a sequence is a basis of the measure space, since the $\sigma$-algebra spanned by all elements of all partitions $\{\eta_n\}$ of such a sequence is dense in the full $\sigma$-algebra.

In the classical encoding, a partition (or a function defining it) is called a {\it generator} if the product of the shifts of the original partition (into the level sets of the function) coincides  $\!\!\mod 0$ with the partition into singletons. In the language of information theory, this means that this encoding loses no information.

Both in the classical and combinatorial cases, the same question arises: does the product of some set of partitions coincide with the partition into singletons? The crucial difference is that in the combinatorial encoding we consider a limit of  {\it finite partitions}, which allows us to use combinatorial tools for solving the problem.

Another difference is in the realization of the quotient by the limiting partition. As we will see, in contrast to the classical case, where  the quotient space is realized as the same symbolic cube $I^{\infty}$ with a new measure and a shift,  in the combinatorial encoding it is realized as the path space of a~locally finite graded graph, or, in other words, as a quasi-stationary (see below) topological chain with finite sets of states, and the quotient endomorphism is realized as a generalized shift, called a transfer.

\smallskip

Of course, the combinatorial encoding in the form described above applies not only to a Bernoulli endomorphism with a continuous set of states, but also to an arbitrary stationary measure in
$I^{\infty}$; moreover, one can start with an automorphism of an abstract Lebesgue space.

\bigskip

Let us state the main problem once again in the most general form.

\begin{problem}[distinguishability problem]
In what cases a sequence of partitions $\{\eta_n\}_n$ separates $\!\!\mod 0$ the points of the space $I^{\infty}$, or is exhausting, or, in other words, when does it solve the distinguishability problem? More formally: in what cases the product $\eta\equiv\bigvee_n \eta_n$ coincides $\!\!\mod 0$ with the partition into singletons (traditionally denoted by $\epsilon$)?
 \end{problem}

The term ``distinguishability'' comes from the fact that the condition introduced above means that {\it almost} any two trajectories $\{\xi_n\},  \{\xi'_n\}\in I^{\infty}$ (from a set of full measure) fall into different elements of the partition $\eta_n$ for sufficiently large~$n$.

Distinguishability is equivalent to the fact that almost every trajectory of the shift can be uniquely recovered from the countable set of elements of the partitions $\eta_n$ containing it. Since the encoding is shift-invariant, it suffices to recover only the first coordinate, all the other coordinates can then be recovered using the shift. The partitions $\eta_n$ are finite, so the positive answer to the distinguishability problem reduces the study of a continuous Bernoulli scheme to that of a countable encoding, i.e., to the study of a sequence of coordinates each taking finitely many values. Recall that a continuous Bernoulli scheme cannot have a finite  generator (in the classical sense), so our construction essentially extends the possibilities of encoding.

\medskip

\subsection{Entropy estimates}

Now we will consider numerical characteristics of the combinatorial encoding in the cases where the distinguishability problem has a positive answer.

Denote by $q_n$ the number of elements in the partition $\eta_n$ and by $r_n(C)$ the number of elements of
the partition  $\eta_{n+1}$ lying in an element $C \in \eta_n$. Recall that
 $r_n(C)\geq 2$, and let $\max_C r_n(C)\equiv r_n$; we have $q_n \geq \prod_1^n  r_k$.

It follows from the definition of the sequence $\{\eta_n\}$ that $S^{-k}\eta_m\prec \eta_{m+k}$. On the other hand, since the entropy $h(S)$ is infinite, none of the partitions $\eta_n$ is a generator, whence
 $\bigvee_n S^{-n}\eta_m \ne \epsilon$. In fact, our aim is to construct a ``diagonal'' refinement of the family of sequences  $ \{S^{-n}\eta_m\}_{n>m}$.

If $\bigvee_n\eta_n=\epsilon$, then every finite partition $\theta$ can be approximated in the entropy metric (see~\cite{Ro}) by a partition that is measurable with respect to $\eta_n$ for sufficiently large $n$; hence, approximating a~sequence $\theta_k$ that approaches a continuous generator and using the invariance of
$\{\eta_n\}_n$, we conclude that
$$
\lim_n \frac{H(\eta_n)}{n}\geq h(S)=\infty.
$$
In particular, if for every $ n$ the partition $\eta_n$ is homogeneous, i.e., all its elements have the same measure, then
      $$\lim_n \frac{\ln q_n}{n}=\infty.$$

Of course, the distinguishability problem will have a positive answer if we allow $q_n$ to grow too rapidly; for instance, it suffices to let $q_n=2^{2^n}$ (a grid with step $1/2^n$ along all $n$ axes of the cube $I^n$); the distinguishability is obvious, since in this case all coordinates are approximated separately. The question is: how slow $q_n$ can grow provided that the distinguishability holds?

On the other hand, we will give a lower  entropy estimate on $\{q_n\}$ which provides a necessary condition for a positive answer to the  distinguishability problem in terms of the growth and can be proved by entropy considerations.

 \begin{proposition}
If for a sequence  $\{\eta_n\}$ the distinguishability problem has a positive answer, i.e.,
$\bigvee_n \eta_n=\epsilon$, then
 $$\lim_n \frac{\ln q_n}{n}=\infty,$$
and the rate of convergence can be arbitrarily small.
 \end{proposition}

But one can give both a direct estimate and a direct construction of a required sequence of partitions.

In the first nontrivial example (the ``Weyl simplices'', see Section~3 below), where $q_n=n!$ and $r_n=n$, the estimate is as follows: $\frac{\ln q_n}{n}=\ln n+O(1)$.
In the second example (the RSK correspondence), the growth is $q_n\sim (n!)^{1/2}$ and $r_n\sim p(n)$, where $p(n)$ is the Euler function (the number of partitions).

  \medskip

In the formal sense, the distinguishability problem in a combinatorial formulation reduces to a~purely computational problem, namely,  checking that  {\it almost all conditional measures on the elements of the partitions $\eta_n$ converge in some metric to $\delta$-measures}, which in other terms is equivalent to  some (nonlinear) law of large numbers. From this viewpoint, the combinatorial encoding is better adapted to the proof of distinguishability than the classical encoding to the proof that some partition is a generator, since in the former case we deal with a limit of finite partitions.

  \section{The main example: encoding  a Bernoulli scheme by Weyl simplices and the triangular compactum $\mathfrak{M}$}

In this paper, we will study the simplest nontrivial example of a combinatorial encoding, namely, the combinatorial encoding of a Bernoulli scheme $(I^{\infty},m^{\infty})$ by Weyl simplices. From the abstract point of view, a special feature of the case under consideration, in terms of a notion introduced above, is that the frame of this encoding is a \emph{homogeneous tree}, in which the number of outgoing edges is the same for all vertices of every level. In this example, i.e., in the ``tree of {\bf i-permutations}\footnote{We introduce the term ``i-permutation'' (coming from ``image of a permutation'') to emphasize
the difference, often neglected, between a permutation as an element of the symmetric group $S_n$ and as an ordering of $n$ objects. If we fix some order on these objects, then an i-permutation is the image of this order under the action of the corresponding permutation.},'' this number is equal to $n$ for vertices of level~$n$, and the number of vertices at this level, as well as the number of paths leading to it, is equal to $n!$.

The main problem is to find the asymptotic behavior of the collection of Weyl chambers (more exactly, Weyl simplices) of type~$A_n$, as well as establish links to related problems.

In fact, this example can be generalized, with the same proofs, to Bernoulli schemes with arbitrary state spaces. However, for definiteness, we will speak mainly of the interval
 $I=[0,1]$ with the Lebesgue measure.

\subsection{The partition of the cube into Weyl simplices}

We will define an increasing shift-invariant sequence of finite cylinder partitions of the cube
$(I^{\infty},m^{\infty})$; namely,  $\eta_n$ is the cylinder partition of~$(I^{\infty},m^{\infty})$ whose base is the partition of the finite-dimensional cube $I^n$ into open Weyl simplices,
by which we mean \emph{the intersections of open Weyl chambers in the Cartan subalgebra with the unit cube}.\footnote{We assume that in ${\mathbb R}^n$ a correspondence is fixed between the Weyl chambers and the i-permutations of the set $\textbf{n}=\{1,2,\dots, n\}$ (in other words, a root system is chosen).}

We will consider the set of full measure in $I^n$ consisting of the vectors with pairwise distinct coordinates. The
\emph{frame} of this sequence of partitions $\{\eta_n\}$ is shown in Fig.~1. This tree can also be called the tree of i-permutations indexing the Weyl simplices (and chambers), or the  \emph{factorial} homogeneous tree.

On the one hand, the set of vertices of level~$n$ consists of all i-permutations of $n$ symbols, which should not be confused with elements of the symmetric group $S_n$. \emph{An edge joins two i-permutations~$a$ and~$b$ if $a$ is obtained from $b$ by removing the element $n$}.

 By definition, the translation (see Section~2.3) associates with every vertex of level~$n$ (for $n>1$) a vertex of level~$n-1$ following the rule according to which the simplex
$\sigma_{x^n}$ of sequences starting from a vector
$x^n=(x_1,x_2, \dots, x_n)$  changes when we remove after application of the shift $S$ the first coordinate $x_1$, that is, pass to the vector
$Sx^n=(x_2,x_3, \dots, x_n)$. Recall that the i-permutation
$(k_1,k_2,\dots, k_n)$ corresponding to the simplex $\sigma_{x^n}$
is given by the formula
    $$k_i=\#\{s \in {\bf n}: x_s<x_i\} i=1,2, \dots n.$$

  \begin{proposition}
The i-permutation $(r_1,r_2,\dots,r_{n-1})$ corresponding to the simplex
$\Sigma_{Sx^n}$ is given by the formula
      $$
 r_i=\begin{cases}
 k_{i+1} &\text{if $k_{i+1}<k_1$},\\
 k_{i+1}-1&\text{if $k_{i+1}>k_1$}.
 \end{cases}
 $$
 \end{proposition}

 The proof immediately follows from definitions. Thus, we have defined a translation, which is a~map from the set of i-permutations of length~$n$ to the set of i-permutations of length~$n-1$. In the next section, where we compute the transfer for this graph, we use this map  and interpret it in a~slightly different way.

Using this rule, we construct the frame corresponding to the tree of Weyl simplices, see Fig.~1.

\begin{figure}
  \includegraphics[width=\linewidth]{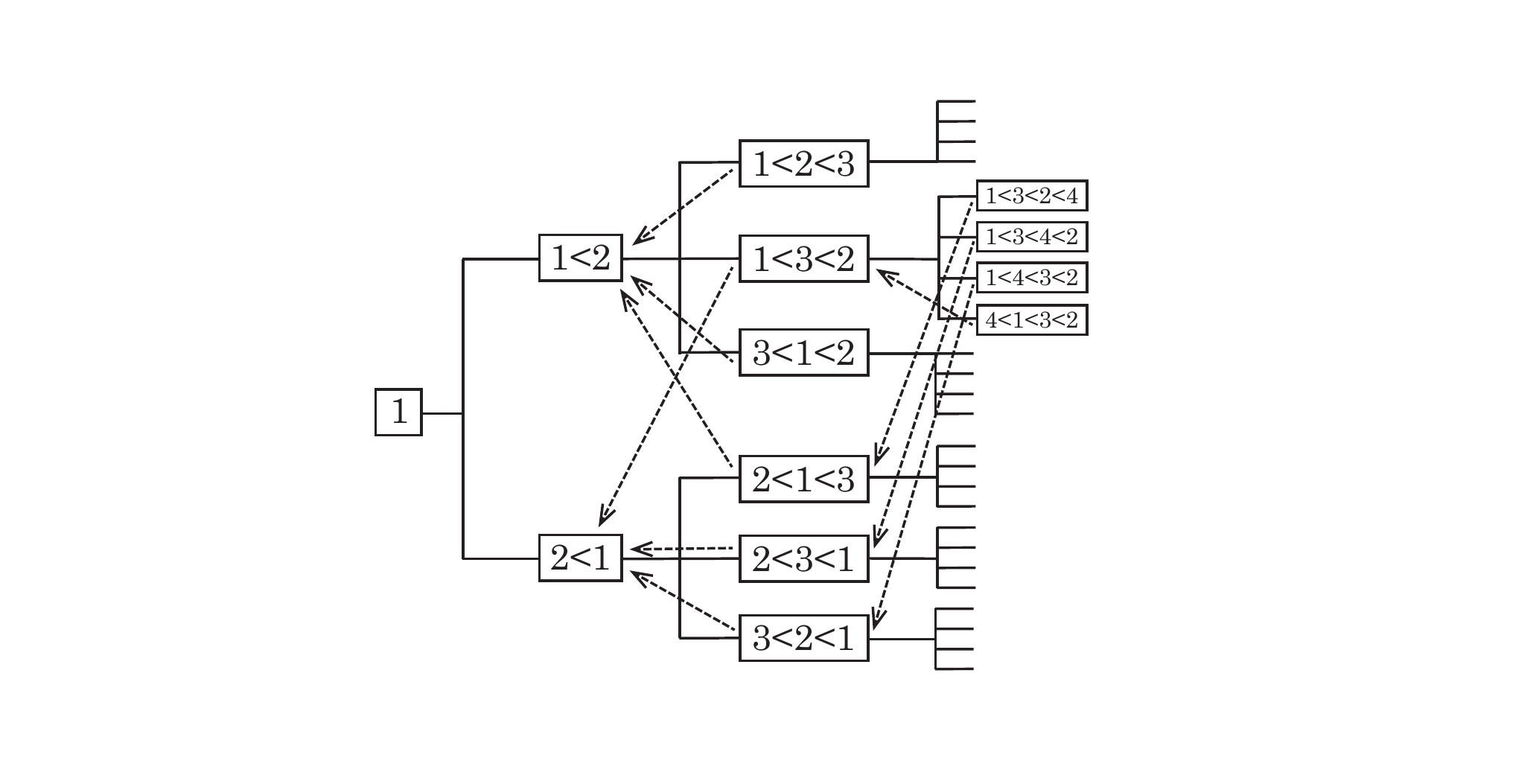}
  \caption{The tree of i-permutations with the translation.}
  \label{fig:Tree}
\end{figure}

The first transition (removing $n$) is called passing to the {\it smaller i-permutation},  and it is natural to say that the translation is passing to the  {\it previous i-permutation}.

It is appropriate at this point to mention the difference between the notions introduced above and the theory of virtual permutations, see~\cite{KOV}: the operation that in~\cite{KOV} and related papers is called passing to the derivative permutation also consists in removing~$n$, but from a permutation rather than an i-permutation
(for example, the derivative permutation of (2413) is (231), while the smaller i-permutation is (213)). That is why, the projective limits with respect to the operations of taking the derivative permutation $S_n\rightarrow S_{n-1}$ and taking the previous i-permutation are different spaces: in the first case, this is the compactum of virtual permutations, while in the second case, we obtain a space whose nature is not quite clear.

Thus, a combinatorial setting for the partition into Weyl simplices is ready. We emphasize the importance of the notion of translation, and hence that of transfer.

An important property of the sequence of partitions of the cube into Weyl simplices follows from the fact that the cube can be represented $\!\!\mod 0$ as the direct product
$$ I^n \sim  (S_n \times \Sigma_n),$$
where $S_n$ is the symmetric group acting on $I^n$ by permutations of coordinates and $\Sigma$ is the standard convex open simplex: $\Sigma=\{{x_1,x_2,\dots, x_n}: 0<x_1<x_2<\dots <x_n<1\}$. The above decomposition is a decomposition of measure spaces, with the normalized Lebesgue measures on
$I^n$ and $\Sigma$ and the uniform measure on $S_n$. Thus, with the partition
$\eta_n$ we can associate its independent complement, which is the partition into the orbits of $S_n$. This direct product decomposition of the cube $I^n$ can be lifted to the infinite-dimensional cube with the Lebesgue measure; namely, the independent complement to the cylinder partition $\eta_n$ is the partition of the cube $I^{\infty}$ into the orbits of $S_n$, which is no longer a cylinder partition.

\subsection{The triangular compactum of paths in the tree of i-permutations}

Consider the ``triangular'' compactum, by which we mean the space of all infinite sequences of positive integers in which the $n$th coordinate takes values in the set
 $\textbf{n}$:
$$\mathfrak{M}=\big\{\{t_n\}_1^{\infty}: t_n\in \textbf{n}=\{1,2,\dots, n\}\big\}.$$
This compactum can be regarded as the set of all paths in the  $\Bbb N$-graded graph $W$ whose $n$th level
is the set $\textbf{n}=\{1,2,\dots, n\}$ and all pairs of vertices of neighboring levels are adjacent. In other words, two neighboring levels in this graph form a complete bipartite graph. The space of all infinite paths of the
 $\Bbb N$-graded graph $W$, denoted by $T(W)$,
is $\mathfrak{M}=T(W)$. But the same compactum is, obviously, the path space of the tree of Weyl simplices.

The compactum $\mathfrak{M}$ resembles the compactum $\mathfrak{S}$  of so-called virtual permutations (see \cite{KOV}), since both are compactifications of the infinite symmetric group, though in different senses; nevertheless, the relation to virtual permutations and representations of the group
 ${\mathfrak{S}}_{\Bbb N}$ is very important and will be discussed below. Recall that the compactum  $\mathfrak{S}$ of virtual permutations was defined as the projective limit of the symmetric groups with respect to the maps
 $$\mathfrak{S}=\lim_{\leftarrow}\{{\mathfrak{S}}_n, p_n\},
 $$
where  $p_n:{\mathfrak{S}}_n \rightarrow {\mathfrak{S}}_{n-1}$ is the operation of deleting the last symbol $n$ from a permutation;  $p_n$ is not a group homomorphism, but it commutes with the right and left actions of  ${\mathfrak{S}}_{n-1}$. Hence, on the projective limit space $\mathfrak{S}$ we have a right and left actions of the infinite symmetric group ${\mathfrak{S}}_{\Bbb N}$. Of course, a virtual permutation is an infinite path in the graph~$W$ of permutations (the graph of Weyl simplices), or a point of the triangular compactum.

One can also define actions of the symmetric group on the triangular compactum, but here we do not discuss this issue. However, an analog of the Haar measure  is defined on $\mathfrak{S}$
and on $\frak M$ in the same way: it is the measure
 $$\mu=\prod_1^{\infty}\mu_n,$$
where $\mu_n$ is the uniform measure on the set $\textbf{n}=\{1,2, \dots, n\}$. This is the unique measure invariant under the left and right actions of the group ${\mathfrak{S}}_{\Bbb N}$.
Besides, there is a one-parameter family of measures $\mu^t$, $t\in \Bbb R_+$, defined as follows:
  $\mu^t(1)=\ldots=\mu^t(n-1)=\frac{1}{t+n-1}$, $\mu^t(n)=\frac{t}{t+n-1}$.
 For $t=1$, we have $\mu^1=\mu$. All measures $\mu^t$ are invariant under the diagonal action.

\subsection{An isomorphism between the cube $I^{\infty}$ and the triangular compactum. The positive answer to the distinguishability problem}

Weyl simplices can be used to define a simple but important isomorphism between the measure spaces
$(I^{\infty}, m^{\infty})$ and $(\mathfrak{M},\mu)$.

Consider the measurable map $J:[0,1]^{\infty}\rightarrow \mathfrak{S}$ from the infinite-dimensional cube
 $I^{\infty}$ to the path space ${T(W)=\mathfrak{M}}$ defined as follows: for $\{x_n\}\in I^{\infty}$,
 $$J(\{x_n\})=\{t_n=t_n(x_1,\dots , x_n)\}_1^{\infty}, \quad\mbox{where}\quad t_n=\#\{i:1\leq i \leq n,\; x_i<x_n\};$$
that is, the $n$th coordinate of the image is equal to the number of coordinates of the preimage with indices at most $n$ whose values are not less than the value of the $n$th coordinate of the preimage.

    \smallskip

We may assume that the map $J$ is defined only on the set of sequences with pairwise distinct coordinates, which has a full measure $m^{\infty}$ in $I^{\infty}$, and is not defined on the remaining set (of zero measure). This is a cylinder map, i.e., the images and preimages of cylinder sets coincide  $\!\!\mod 0$ with cylinders in the corresponding spaces.

Consider in more detail the map $J$ for a finite-dimensional cube. Obviously, the $J$-preimage of a~point of the space $\prod_{i=1}^n\textbf{i}$ is an open Weyl simplex, i.e., the intersection of an open Weyl chamber with the unit cube: this is the open simplex of all vectors of the unit cube with a fixed collection of pairwise inequalities between their (distinct) coordinates. Thus, the finite-dimensional level sets of the map $J$ divide ($\bmod\,0$)
 the unit cube into the Weyl simplices.

\begin{theorem}
The map $J$ is a $\bmod\,0$ measure-preserving isomorphism between the spaces
 $(I^{\infty},m^{\infty}) $ (infinite-dimensional cube) and $(\mathfrak{M},\mu^{\infty})$ (triangular compactum). It sends the sequence $\{\eta_n\}$ of partitions of the space
$(I^{\infty},m^{\infty})$  defined above to the sequence $\tau_n$ of complete cylinder partitions of the compactum $\mathfrak{M}$.
\end{theorem}

To make it clear, $\tau_n$ is the partition of $\mathfrak{M}$ into the classes of sequences of positive numbers in which the first $n$ coordinates coincide, $\prod_n \tau_n=\epsilon$.

\medskip

One can see from the structure of the map~$J$ that the values of the measures $\mu^{\infty}$ and $m^{\infty}$ on cylinder sets agree; indeed, both the measure of a finite-dimensional Weyl simplex of order $n$ and the measure of a point in a finite $n$-fragment of $\mathfrak{M}$ are equal to $\frac{1}{n!}$, which implies that the measures of preimages and images coincide, i.e., the map $J$ is defined almost everywhere in  $(I^{\infty},m^{\infty}) $, is  $\bmod\, 0$ surjective, and preserves the measure.

It remains to prove that $J$ is $\bmod\, 0$ injective, i.e., separates almost all points of the preimage
 $I^{\infty}$. This simple but remarkable statement is worth highlighting.

\begin{lemma}
The limiting partition $\eta=\lim_n \eta_n$  of the infinite-dimensional cube $I^{\infty}$ (the limit of the partitions into open Weyl simplices) coincides
 $\bmod\, 0$ (with respect to the Lebesgue measure) with the partition into singletons. In other words, the distinguishability problem for the partition into Weyl simplices has a positive answer. Therefore, the map $J$ is an isomorphism of measure spaces. In more detail, there exists a set of full Lebesgue measure in $I^{\infty}$ such that for any two points~$\{x_n\}$ and~$\{x'_n\}$ of this set there exist indices
$i$ and $j$ for which the corresponding coordinates satisfy the opposite inequalities:
$$x_i > x_j, \quad \mbox{but}\quad x'_i <x'_j.$$
\end{lemma}

In a somewhat paradoxical form, the lemma can be stated as follows: almost every (with respect to the Lebesgue measure) infinite sequence of points from the interval $[0,1]$ can be uniquely recovered from the list of pairwise inequalities between its coordinates.

\smallskip

Or, even more paradoxically: almost every infinite-dimensional Weyl simplex consists of a single point.\footnote{If, instead of the cubes $I^n,I^{\infty}$  we consider the spaces ${\Bbb R}^n, {\Bbb R}^{\infty}$ with the standard infinite-dimensional Gaussian measures, then our statement looks as follows: almost every infinite-dimensional Weyl chamber consists of a single ray.}

\begin{proof}
Assume that two sequences have the same inequalities for all pairs of coordinates but differ in at least one coordinate: $x_n=\alpha < \beta=x'_n$. Since the coordinates are independent, it follows that with probability~$1$ there exists a number $N$ such that
 $$x_N\in (\alpha,\beta), \quad\quad x'_N\in (\alpha,\beta),$$ which implies a contradiction:
      $\alpha=x_n<x_N$, but  $\beta=x'_n> x'_N$.
\end{proof}

In fact, the proof uses not the independence of coordinates, but the equidistribution of almost every sequence  $\{x_n\}_n$, which follows from the pointwise ergodic theorem; hence, we can replace the Lebesgue measure with any measure in $I^{\infty}$ for which the coordinates are equidistributed.

Let us sketch another, more conceptual argument, which is applicable in a much more general situation.

To prove the distinguishability of Bernoulli realizations (trajectories) for a sequence of partitions means to prove the following: for every measurable function $f$ (it suffices to consider only cylinder functions), its average over the conditional measure on an element $C$ of the partition
 $\eta_n$ approaches~$f$  for large $n$. This, in turn, is equivalent to the fact that the projection of the conditional measure on the element $C$ to the finite-dimensional simplex whose cylinder hull contains $C$ is close to the $\delta$-measure at the barycenter of this simplex. But this is true by the equidistribution of almost all trajectories and their fragments. Indeed, the equidistribution implies that the projections of long fragments of the sequence concentrate near the barycenter, since the interval between the coordinates should be filled uniformly. These considerations suffice to recover this proof of the theorem.

 In contrast to the proof given above, the last argument does not use specific features of the partition into Weyl simplices; it is universal, since everything reduces to equidistribution or, more generally, to a law of large numbers. Elsewhere, we will apply the same argument to obtain a new proof of the Romik--Sniady theorem  \cite{RS,Sn}  (which, in our terms, is a theorem on the distinguishability of the encoding via $Q$-tableaux in the RSK correspondence).

\subsection{Transfer for the tree of i-permutations and the triangular compactum}

Now we must write the image of the shift $S$ under the isomorphism $J$ in terms of sequences
 $\{t_n\}_n \in \mathfrak{M}$, i.e., find the corresponding transformation of the triangular compactum
 $\mathfrak{M}$:
    $$\Lambda= JSJ^{-1}: \mathfrak{M} \rightarrow \mathfrak{M}.$$
This is exactly the transfer of the triangular compactum regarded as the path space of the graph    $W$. To find it, we will use  the formula we have obtained for the translation.

Denote $x^n=\{x_i\}_{i=1}^n$, and let $d_n(x^n)$ be the number of coordinates in $x^n$ that are less than $x_1$. Clearly, each $d_{n+1}$ is equal either to $d_n+1$ (if $x_{n+1}<x_1$), or to $d_n$ (if $x_{n+1}>x_1$). It is convenient to use the following terminology.

For each finite fragment of a path $t=\{t_n\}$ in the image (i.e., in the compactum $\mathfrak{M}$), we define  \textbf{marked positions} by induction as follows.
The first position $t_1=1$ is marked by definition. Assume that the number of marked positions among the first $n$ coordinates is equal to $d_n(t)$; then the position~$t_{n+1}$ is marked if and only if
$t_{n+1}\leq d_n$, i.e.,  $t_{n+1}$ does not exceed the number of previously marked positions.

\begin{theorem}
The formula for $\Lambda= JSJ^{-1}$ is as follows:
 $\Lambda(\{t_n\})=\{t'_n\}$
 where
 $$
 t'_n=\begin{cases}
 t_{n+1} &\text{if $t_{n+1}$ is marked},\\
 t_{n+1}-1&\text{if $t_{n+1}$ is not marked.}
 \end{cases}
 $$
\end{theorem}

In other words, the $n$th coordinate   $t'_n$ of the image either coincides with the $(n+1)$th coordinate~$t_{n+1}$ of the preimage, or is less  by one, depending on whether the number of coordinates less than the first one increases by~$1$ when we add the $(n+1)$th coordinate in the preimage.

In short, the transfer sends a virtual permutation to a new virtual permutation in which the $n$th position is occupied either by the number that occupied the $(n+1)$th position in the original permutation, or by this number decreased by~$1$, depending on its value.

Besides, the following relations hold:
  $$t_{n+1}-t'_n =1-(d_{n+1}(t)-d_n(t)).$$
The proof immediately follows from the previous formulas and considerations.

Let us turn to the inversion formula.  The formulas for $d_n$ directly imply the following theorem.

\begin{theorem}
For almost all trajectories $\{x_n\}_n\in I^{\infty}$ with respect to the measure $m^{\infty}$,
$$\lim_n \frac{d_n}{n}=x_1.$$
 \end{theorem}

In the same way we can find the other coordinates $x_n$, $n>1$. This and other similar formulas can be regarded as inversion formulas for the isomorphism $J$.

Thus, we have completely described an isomorphism $J$ between the triples
$$(I^{\infty},m^{\infty}, S)\quad \mbox{and}\quad (\mathfrak{S},\mu,\Lambda).$$

The action of the operator $\Lambda$ on the space of virtual permutations with the Haar measure (which is isomorphic to a Bernoulli action) is of interest. One can prove the Bernoulli property for this operator directly (i.e., without using the isomorphism $J$), and also present a Bernoulli generator, which is ``expelled to infinity,'' as shown by the formula from Theorem~3. Here we see a remote analogy with the boson-fermion correspondence in a combinatorial version. The operator $\Lambda$ is exactly the \emph{transfer} defined in the previous section.

The map $J$ also establishes an above-mentioned isomorphism between the sequence of partitions~$\{\eta_n\}_n$ of the space $I^{\infty}$  and the sequence of partitions $\{\tau_n\}_n$ of the space $\mathfrak{S}$.

Let us summarize our considerations.

\begin{theorem}
The one-sided Bernoulli shift $S$ on the space $(I^{\infty}, m^{\infty})$ is metrically isomorphic to the transfer
$\Lambda$ defined on the triangular compactum (the compactum of virtual permutations) with the Haar measure. An isomorphism $J$ is established by the encoding of the cube with the system of Weyl simplices, which form an exhausting increasing sequence of finite partitions.

Another model of the triangular compactum is the path space of the infinite homogeneous tree of permutations with the transfer and the uniform central measure.
\end{theorem}

A nontrivial property of the isomorphism $J$ is that it reverses the direction of time: the first coordinate (as well as other cylinder functions  on the cube) is mapped by this isomorphism to an infinitely remote limiting function on the triangular compactum, which is not a cylinder function; though some cylinder functions remain cylinder functions, but their order increases. The continuity of the original system ``escapes'' to infinity.

\section{Comments and remarks}

Several comments are in order related to the above example.

\subsection{More on the distinguishability problem}

Using the combinatorial encoding, we have obtained an example of a realization of a Bernoulli endomorphism as a shift in the path space of a graph or, in other words, as a shift in the space of a~nonstationary chain. This realization is a special case of the notion of transfer.

Another special, but much more complicated case of transfer is the Sch\"utzenberger transformation. It arises when one considers a covering by a Bernoulli scheme of the space of infinite standard Young tableaux with the Plancherel measure, which was (in the general case of ergodic central measures) suggested in~\cite{KV}. A deep analysis and a proof of the fact that this covering is an isomorphism (rather than only a homomorphism) between the Bernoulli shift and the Sch\"utzenberger transformation, i.e., in our terms, a proof of distinguishability, was recently given by
 D.~Romik and P.~Sniady~\cite{RS,Sn}. Their analysis is based on the study of the Sch\"utzenberger transformation from the viewpoint of what I~have called the nerve of a tableau using limit shape techniques, which allow one to obtain a complicated inversion formula. We will consider the distinguishability problem in the general setting, in particular, for the RSK correspondence, in another paper.

\subsection{The list of isomorphisms}

Let us enumerate isomorphic spaces with an invariant measure and a transfer:
\begin{itemize}
\item the path space of the factorial tree;
\item the path space of the graph~$W$;
\item the triangular compactum with the Lebesgue measure;
\item the space $(I^{\infty}, m^{\infty})$ with the Bernoulli endomorphism.
\end{itemize}
This list can be extended by other graphs and their path spaces. For instance, by the Young graph and its path space (i.e., the space of infinite standard Young tableaux) with the Plancherel measure.

\subsection{The relation to matrix distributions}

Consider the space $M_{\Bbb N}^{symm}(\pm 1)$ of all infinite symmetric matrices with zeros on the principal diagonal, and the product measure $\sigma$ on it such that all coordinates
 $\epsilon_{i,j}$, $i>j$, are independent and have the distribution $(1/2,1/2)$.

\begin{theorem}
Let a map $T:I^{\infty}\rightarrow M_{\Bbb N}^{symm}(\pm 1)$ be defined by the formula
$$\epsilon_{i,j}
 =\begin{cases}
 +1&\text{if $x_i>x_j$},\\
 -1&\text{if $x_i<x_j$}.
 \end{cases}
 $$
 It is a metric isomorphism between the above spaces endowed with the measures $m^{\infty}$ and $\sigma$, respectively.
\end{theorem}

The fact that $T$ is an isomorphism of measure spaces follows from the lemma proved above saying that the limit of the partitions $\eta_n$ coincides $\bmod\, 0$ with the partition into singletons: $\lim_n \eta_n=\epsilon$.

For completeness, observe that the measure $\sigma$ on the space of matrices is the matrix distribution (in the sense of~ \cite{FA}) of the following measurable function on the unit square endowed with the Lebesgue measure:
$\phi(x,y)=\operatorname{sgn} (x-y)$. It is invariant under the simultaneous action of the infinite symmetric group on the rows and columns of matrices.

\medskip

In fact, above we have established a not quite obvious isomorphism between the space of such matrices endowed with the measure $\sigma$ and the space of virtual permutations $\mathfrak{S}$ endowed with the Lebesgue measure $\mu$. The measure~$\sigma$ has interesting properties, it has also appeared in other contexts (universality), see~\cite{GW}.

\subsection{Other compactifications of symmetric groups}

The notion of translation of i-permutations  suggests the idea of considering the projective limit of the symmetric groups with respect to this operation. We will assume that it applies to permutations
rather than i-permutations, i.e., under the identification of i-permutations with permutations, maps
 $S_n$ onto $S_{n-1}$. This map resembles the operation of taking the derivative permutation, but is more complicated. Namely, we remove the first element (the image of~$1$), and then decrease the coordinates that are greater than the removed one by $1$, leaving the other coordinates unchanged. The projective limit of the symmetric groups $S_n$ with respect to these operations should be regarded as a completion of the infinite symmetric group.

It seems that apart from this operation and that of taking the derivative permutation (see \cite{KOV}), there is a whole series of reasonable operations and the corresponding projective limits. For example, one can delete~$1$ in the same way as $n$ is deleted when taking the derivative permutation; this also results in an interesting object, a new extension of the infinite symmetric group different from the space of virtual permutations. Here is a more detailed description of the map that deletes~$1$ from a~permutation and decreases all numbers by~$1$:
$$
S_n \ni g=\begin{array}{ccccccc}
  1   &  2  &  3  & ... & r & ... & n \\
  r_1 & r_2 & r_3 & ... & 1 & ... & r_n
\end{array}
\rightarrow Tg=
\begin{array}{ccccccc}
    1   &   2   & ... &  r-1  &...& n-1  \\
  r_2-1 & r_3-1 & ... & r_1-1 &...&r_n-1
\end{array}\in S_{n-1}.
$$

Obviously, this map is defined as a map (but not a homomorphism) from the group
${\frak S}_{\infty}$ to itself such that every element has countably many preimages. It seems that this ``nonhomomorphism'' of the group ${\frak S}_{\infty}$ to itself has every right to exist.

\subsection{Relation to $C^*$-algebras}

From the viewpoint of the general theory, the example considered above corresponds to the Glimm
$C^*$-algebra $\bigotimes_{n=1}^{\infty}
M_n(\Bbb C)$. Indeed, the spectrum of a finite-dimensional diagonal subalgebra of this graded algebra is exactly a finite fragment of the triangular compactum, and the whole compactum is the spectrum of the Gelfand--Tsetlin algebra for this $C^*$-algebra. The unique central measure on  the triangular compactum is the Haar measure.

\subsection{Representations of the infinite symmetric group}

The established isomorphism between the space of virtual permutations and the infinite-dimensional cube allows one to obtain new models of representations of the infinite symmetric group. It is more important that the result obtained above establishes an asymptotic isomorphism between the  $l^2$-space $l^2(\mathfrak{S}_{\infty})$ of the infinite symmetric group  $\mathfrak{S}_{\infty}$ and a space of $L^2$-functions on the Cartan subalgebra of a Lie algebra of type~$A_n$.

For every $n$, we define an embedding of the space $l^2(S_n)$ into the space of $L^2$-functions on the distinguished Cartan subalgebra
${\frak H}_n$ of the Lie group $GL_n(\Bbb R)$ with respect to the Lebesgue measure on the unit cube
$I_n \subset {\frak H}_n$ in the chosen coordinates.\footnote{We can just as well consider the $L^2$-space
with respect to the standard Gaussian measure on ${\frak H}_n$ corresponding to the Killing form.}

An embedding $\varrho_n: l^2(S_n) \rightarrow L^2_m(X_n)$ is defined by the formula
  $$\varrho_n(\delta_w)=\textbf{1}_{\Delta_w},$$
where $\delta_w \in C(S_n)\subset l^2(S_n)$ is the $\delta$-function at an element $w\in S_n$ and  ${1}_{\Delta_w}$ is the characteristic function of the Weyl simplex $\Delta_w$. This correspondence can be extended by linearity to the whole group algebra of~$S_n$, and we obtain a map from this algebra to the space of functions on the Cartan subalgebra, which will be denoted by the same symbol~$\varrho_n$.

\begin{proposition}
The embedding $\varrho_n$ is an isometry with respect to the norms of the spaces $l^2(S_n)$ and~$L^2_m(X_n)$.
\end{proposition}

Note that for every $n$ the image of the space $l^2(S_n)$ under the isometry $\varrho_n$ is a proper finite-dimensional subspace in $L^2_m(X_n)$.

It is natural that the $l^2$-space on the Weyl group can be embedded into the space of functions on the Cartan algebra, but it is by no means obvious that in the limit the image coincides with the whole space.

It turns out that this is true by the positive answer to the distinguishability problem: in the limit we obtain an isometry. In more detail, we assume that $l^2(S_n)$ is, in a natural sense, a subspace in~$l^2(\mathfrak{S}_{\infty})$, and the spaces of $L^2$-functions on the Cartan subalgebras also constitute an inductive family, since we have embeddings of groups. Hence we can consider the limit of the isometries $\varrho_n$.
\begin{corollary}

The limit of the isometries
$ \varrho_n:l^2(S_n)\rightarrow
L^2_m(I_n)$ exists and is an isometry between the spaces $l^2(\mathfrak{S}_{\infty})$  and $L^2(I^{\infty})$. It follows that the regular representation of the group
$\mathfrak{S}_{\infty}$ can be realized, in a natural way, in the $L^2$-space over a Bernoulli scheme as the limit of actions of finite Weyl groups on Cartan subalgebras.
\end{corollary}

This realization and its generalizations should be studied further.  As we have already mentioned, such a realization of a Bernoulli endomorphism exists for a whole series of graphs, and each of them corresponds to a realization of the representation of the infinite symmetric group which in~\cite{V} was called the basic representation.

\subsection{An analog of the partition into Weyl simplices for an arbitrary Bernoulli scheme}

We conclude this section by showing that the method of encoding via Weyl simplices can be used to encode an arbitrary sequence of independent variables, and the linear order on the interval can be replaced by an arbitrary measurable ordering on a set of full measure.

Consider an arbitrary Lebesgue space $(I,m)$
and divide its square $I \times I$ into a measurable set $A$ and its complement $\bar A$ satisfying the following condition: for every pair of points $u,v \in I$ there exists a set $X_{u,v} \subset I\times I$ of positive measure such that for every $w\in X$  one of the pairs $(w,u), (w,v)$ lies in~$A$ and the other one lies in $\bar A$.

 \begin{theorem}
 Consider a Bernoulli scheme with the state space
 $(I,m)$ and assume that we have chosen a set $A\subset I\times I$ of measure~$1/2$ satisfying the condition stated above. Construct an increasing invariant sequence of finite measurable partitions $\{\eta_n\}$ as follows: an element $C$ of the partition
 $\eta_n$ is a cylinder with base $I^n$ consisting of all collections $(x_1,x_2,\dots,x_n)$, $x_i\in I$, for which every pair $(x_i,x_j)\in I\times I$ for $i,j=1,2,\dots, n$
lies in a fixed one of the sets
 $A$ or $\bar A$.
Then

 {\rm1.} The sequence $\{\eta_n\}$ has the same frame as the sequence of Weyl simplices and is metrically isomorphic to it. In particular,

{\rm2.} The distinguishability problem for this sequence of partitions $\{\eta_n\}$ has a positive answer.
 \end{theorem}

The proof is exactly the same as for Weyl simplices. Thus, from the metric point of view, the linear order on the interval has no specific features.

\subsection{Distinguishability of random matrices}
The statement  of the distinguishability problem and its analysis suggested above can be used in a~large number of problems similar to those considered in the paper. Let us give an example from the theory of random Gaussian matrices (GOE). Consider the space $M_{\Bbb N}(R)$ of all infinite symmetric real matrices endowed with the standard Gaussian measure  $\mu$. Here we consider the measure space $(M_{\Bbb N}(R),\mu)$ instead of  $(I^{\infty}, m^{\infty})$. Define a (no longer increasing) sequence of (no longer finite)  measurable partitions~$\{\eta_n\}_n$ of this space: an element of~$\eta_n$ is the set of all matrices for which the  $n\times n$ principal minors  have the same spectrum. The parameter of~$C$ is an $n$-vector $\lambda^n(C) $ of eigenvalues of the corresponding submatrix.

We obtain the
{\it distinguishability problem for the sequence $\{\eta_n\}_n$: does there exist a set of matrices of full measure $\mu$ that can be uniquely recovered from the collection of parameters
$\{\lambda^n(C)\}_{C,n}$}? Cf.\ Lemma~1 on the distinguishability of points
$\{\xi_n\}\in I^{\infty}$ by the collection of pairwise inequalities.

\section{Graded graphs, transfer, and quasi-stationary processes}

\subsection{Why not only trees}

The combinatorial encoding of a sequence of random variables involves the study of increasing invariant sequences of finite partitions, and the combinatorial counterpart of the problem reduces to the study of trees with a transfer. However, in some natural settings of the distinguishability problem, the original object is not necessarily an increasing sequence of partitions.

Assume that we have a sequence of finite partitions $\zeta_n$  of a Lebesgue space $(X,m)$. If we regard the set of all elements of all partitions as the set of vertices of a graph to be constructed, and assume that an element $C$ of $\zeta_n$ does not necessarily lie in one element of $\zeta_{n-1}$, but may intersect several elements $D$ of this partition, then we obtain a graph with edges connecting elements of neighboring partitions (i.e., vertices of neighboring levels) $C$ and $D$ if and only if the measure of their intersection is positive. Thus, we obtain a graded graph in which for each vertex there is a
probability vector on the set of edges entering this vertex  (i.e., an equipped graph in the terminology of~\cite{UMN}). The tree defined in Section~2 as the frame of an increasing sequence is a special case of this construction.

If we pass from the partitions $\zeta_n$ to their products $\bigvee_{k=1}^n\zeta_k=\eta_n$, then we obtain an increasing sequence of partitions $\{\eta_n\}_n$, which reduces the problem to the case of a tree considered above; moreover, in terms of the graph constructed above, this tree is nothing else than the tree of finite paths in the original graph leading from the initial vertex to all the other vertices, and the graph itself is a natural quotient of this tree.

However, this reduction does not at all mean that there is no benefit from considering the graph itself, which can be illustrated by the example of encoding a sequence of independent variables by the Young graph and the RSK correspondence, as in~\cite{KV,Sn,RS}. In other words, the distinguishability problem can be analyzed in terms of the graph itself without passing to the tree of paths.

\subsection{Definition of transfer for a graded graph}

Let us sketch the main idea of applying the theory of graded graphs to ergodic problems of the type under consideration. The key role here is played by the notion of transfer in the path space of a graded graph which we introduce below.

\medskip
Consider an arbitrary Bratteli diagram, i.e., an $\Bbb N$-graded locally finite graph (or even a multigraph)~$\Gamma$. An infinite tree is an example of such a graph. A path in~$\Gamma$ is an infinite maximal sequence of edges in which the beginning of each edge coincides with the end of the previous edge. Denote the space of all paths by  $T(\Gamma)$; this is a Cantor-like compactum in the inverse limit topology. A~transformation
  $$\Lambda: T(\Gamma)\rightarrow T(\Gamma)$$
  is called a \emph{transfer} if it is continuous, decreases the level of each edge by~$1$, and satisfies the following locality (Markov) property:  for every $n$, the rule according to which an edge between levels~$n+1$ and~$n+2$ in a path is mapped to an edge between levels~$n$ and~$n+1$ in the image of this path depends only on the fragment of the path between levels~$n$ and~$n+2$.  This means that a transfer is determined by a set of local rules for the translation of an edge to the previous level.

For stationary graphs, in which the sets of vertices of all levels (except the first one) are isomorphic and these isomorphisms are fixed, the translation rule depends on nothing: an edge connecting vertices~$a$ and $b$ of levels~$n+1$ and~$n+2$ goes to the edge connecting the vertices~$a'$ and~$b'$ of levels~$n$ and~$n+1$ identified with the vertices~$a$ and~$b$, respectively. In this case, the transfer is an ordinary shift.

For trees regarded as graded graphs, the definition of transfer coincides with that from Section~2.3.

\begin{definition}
A graded graph is said to be quasi-stationary if a ``transfer'' operation,  corresponding to some translation of edges in the sense described above,  is defined on its path space. The path space of a~quasi-stationary graph, regarded as a topological Markov compactum, will be called a quasi-stationary Markov compactum.
\end{definition}

Thus, we have described a new type of realizations of automorphisms and endomorphisms with infinite entropy as transfers on quasi-stationary Markov compacta.

According to our definition, a transfer  is a shift of sequences of edges, and not of sequences of vertices as in the stationary case. Thus, this notion opens new possibilities for realizations of transformations.

\begin{figure}
  \includegraphics[width=\linewidth]{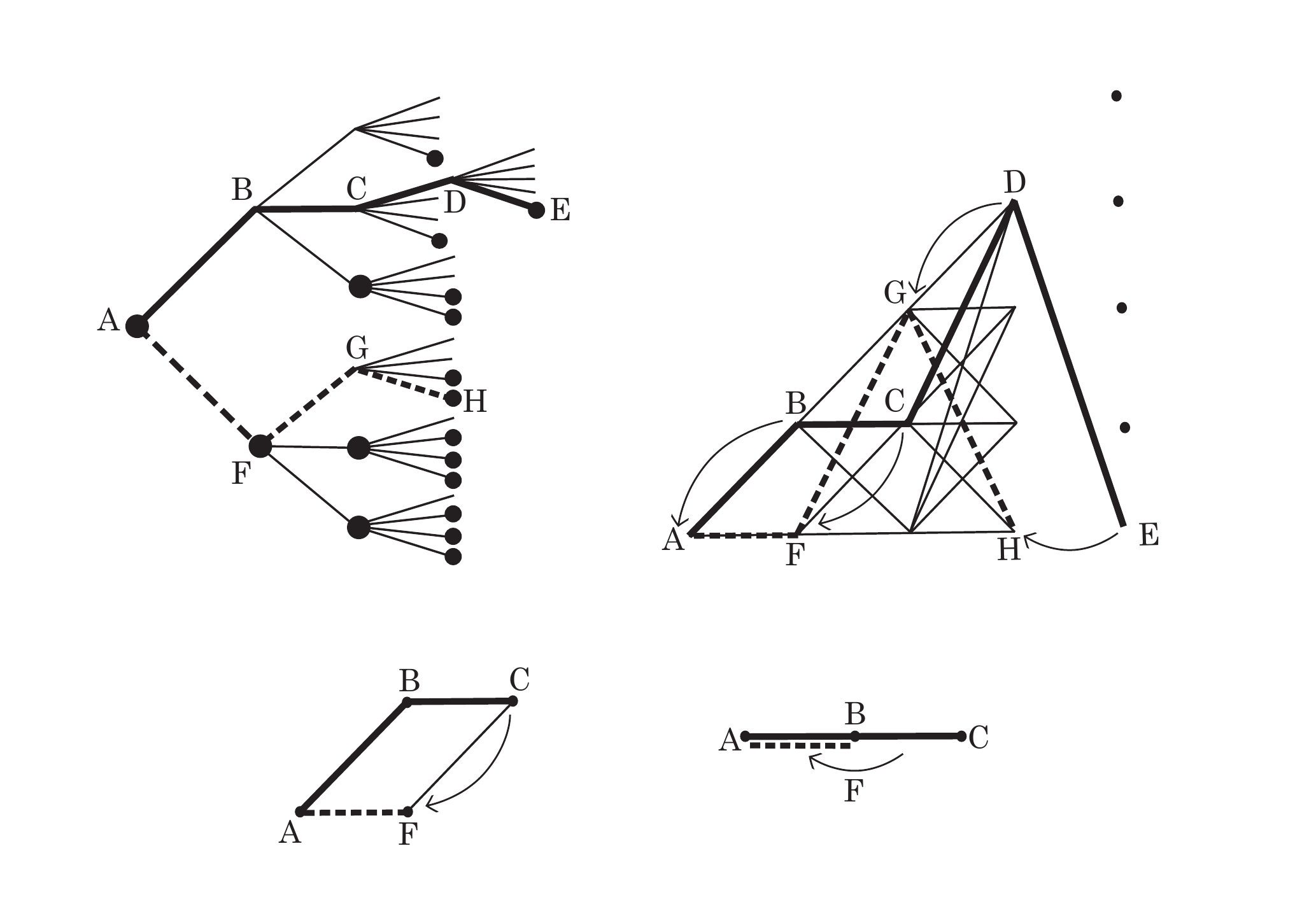}
  \caption{The dashed path is the transfer of the bold one.}
  \label{fig:transfer}
\end{figure}

A transfer defines an additional structure on the graph and, in general, is not uniquely determined by the graph itself, though in some cases there exists a distinguished transfer.
For example, in the important special case of graphs in which every $2$-interval contains either one or two intermediate vertices (this is the case for Hasse diagrams of arbitrary distributive lattices, in particular, for the Young graph and other examples), a translation of edges between adjacent vertices is determined in a~natural way by the very structure of the graph, see Fig.~2.

\begin{proposition}
For the Young graph, a transfer on the path space (i.e., on the space of infinite standard Young tableaux) is defined automatically in the sense described above, since the Young graph is the Hasse diagram of the distributive lattice of finite ideals of the lattice ${\Bbb Z}_2$.  In this case, it coincides with the well-known Sch\"utzenberger transformation, which is a special case of transfer.
\end{proposition}

The proof follows from a detailed analysis of the definition of transfer.
(For the Sch\"utzenberger transformation, see \cite{St,KV,Sn,RS}.)

If  a transfer is defined on the path space of a graded graph, then this space should be regarded as a~nonstationary (or {\it quasi-stationary}) Markov chain, meaning that  the transfer is an analog of the shift. If we have a central measure on the path space that is invariant under the transfer, then we obtain a~quasi-stationary Markov chain with an invariant measure.  Hence the theory of transfer becomes part of ergodic theory, as a nonconventional realization of measure-preserving transformations. A~more detailed exposition of the theory of transfer will be presented elsewhere.

\medskip

The author is grateful to P.~P.~Nikitin for preparing the figures, P.~B.~Zatitskii for reviewing the literature
and G.M.Zukerman for creative attitude to the text.

Translated by N.V.Tsilevich.

\end{document}